\documentclass[a4paper,12pt]{amsart}
\usepackage{amsmath,amsfonts,amssymb}
\usepackage{amsthm}
\usepackage[usenames]{color}
\usepackage{cite}
\usepackage{color}
\usepackage{hyperref}
\usepackage{tikz}

\newtheorem{theorem}{Theorem}[section]

\newtheorem{proposition}[theorem]{Proposition}
\newtheorem{lemma}[theorem]{Lemma}

\newtheorem{remark}[theorem]{Remark}

\renewcommand{\proof}{{\noindent \bf Proof:\ }}

\numberwithin{equation}{section}

\title[Pullback attractor for a non local non-autonomous  evolution equation]{Pullback attractor for a non local non-autonomous  evolution equation in an unbounded domain}

\author[F. D. M. Bezerra]{Flank D. M. Bezerra$^1$}\thanks{$^1$Research partially
supported by FAPESP \# 2011/04166-5, Brazil}
\address[F. D. M. Bezerra]{Universidade Federal da Para\'{\i}ba, Departamento de Matem\'atica, 58051-900, Jo\~{a}o Pessoa, PB, Brazil.}
\email{flank@mat.ufpb.br}

\author[M. S. Pereira]{Miriam da S. Pereira$^2$}
\address[M. S. Pereira]{Universidade Federal da Para\'{\i}ba, Departamento de Matem\'atica, 58051-900, Jo\~{a}o Pessoa, PB, Brazil.}
\email{miriam@mat.ufpb.br}

\author[S. H. da Silva]{Severino H. da Silva$^3$}\thanks{$^3$Research partially
supported by CAPES/CNPq, Brazil}
\address[S. H. da Silva]{Universidade Federal de Campina Grande, Unidade Acad\^emica de Matem\'atica e Estat\'istica,  58429-900, Campina Grande, PB, Brazil.}
\email{horacio@dme.ufcg.edu.br}

\date{\today}

\begin{document}

\maketitle

\begin{abstract}
In this work we consider the non local evolution equation with time-dependent terms which arises in models of phase separation in $\mathbb{R}^N$
\[
\partial_t u=- u + g \left(\beta(J*u)   +\beta h(t,u)\right)
\]
under some restrictions on $h$, growth restrictions on the nonlinear term $g$ and $\beta>1$. We prove, under suitable assumptions, existence, regularity and upper-semicontinuity of pullback attractors with respect to functional parameter $h(t)$ in some weighted spaces.

\vskip .1 in \noindent {\it Mathematical Subject Classification 2010:} 35B40, 35B41, 37B55.\\
\noindent{\it Keywords}: pullback attractors; upper semicontinuity; nonlocal evolution equation; weighted Sobolev spaces.
\end{abstract}

\section{Introduction}

The continuum limit of one dimensional Ising spin systems with Glauber dynamics and Kac
potentials gives rise, see \cite{AMRT_1, Masi_1} and \cite{Masi_3}, to the non local evolution equation
\begin{equation}\label{ref1}
\partial_t u(x,t)=- u(x,t) + \tanh \left(\beta(J*u)(x,t)   + \beta h\right),
\end{equation}
where $u(t,x)$ represents the magnetization density in $x\in\mathbb{R}$ at time $t\in[\tau,+\infty)$; $(J*u)(x)=\int_{\mathbb{R}} J(x-y)u(y,t)dy$; $\beta>0$ the inverse
temperature of the Ising system; $J\in\mathcal{C}^1(\mathbb{R})$ a non-negative even function which gives the strength of the spin-spin interaction; $h$ an constant external magnetic field.

Non local equations like (\ref{ref1}) are well studied in physics of phase separation and interface dynamics, see for instance  \cite{KUH1}, \cite{KUH2} and \cite{KUH3}, but also in many other fields as biology, population dynamics, propagation of diseases, see for instance \cite{DK} and \cite{FM}.

In the theory of dynamical systems in infinite dimensional spaces this equation also have been studied widely in different contexts. For instance, in \cite{Masi_1} the authors study the case of an external magnetic field $h$, characterizing the travelling-front solutions of (\ref{1.1}) for small values of $h$, and proving that their shape is globally stable. Here, we consider $h^*>0$ a constant such that the equation
\begin{equation}\label{Def-h}
s=\tanh(\beta s+\beta h(\cdot,s)),\ \ 0\leqslant h<h^*
\end{equation}
has three and only three different roots for $s\in[-1,1]$.

In \cite{Be}, \cite{Silva} and \cite{severino2} the analysis of the asymptotic behavior of solution for non local diffusion equations is performed under the point of view of the theory of compact global attractors. In \cite{Pereira} the author prove the existence of a global attractor for this equation in some weighted spaces in the context of one dimensional, the author also studies the characterization of such attractors with the existence of nonhomogeneous equilibria.

In this paper, we are concerned with the study of the asymptotic behavior of solutions to initial value problems associated with non-autonomous equations motivated by (\ref{ref1}), see (\ref{1.1}) below. These equations can be seen as a non-autonomous ODEs in Banach spaces, and therefore the properties of (local) existence  and uniqueness follow from standard results of the classical theory. Our interest in the problem comes from the fact that the solutions of these problems shares with the solutions of semilinear parabolic (or hyperbolic) non-autonomous problems interesting qualitative properties, such as the existence of smooth pullback attractors. However, the investigation of qualitative properties of the evolution process given by these equations is a much harder topic.

The essential difference between the results here and the we mentioned before is that our goal is to prove, under some hypotheses about the function $g$ on $\mathbb{R}$, the existence, regularity and upper semicontinuity of pullback attractors for the non local model (\ref{1.1}). With this, our work completes the study of equation with nonlocal terms.

More specifically, in the Banach space $L^2(\mathbb{R}^N)$, we consider the following  non local non-autonomous evolution equation
\begin{equation}\label{1.1}
\partial_t u(t,x)=
- u(t,x) + g \left(\beta(J*u)(t,x))   +\beta h(t,u(t,x))\right)\quad\mbox{and}\quad u(\tau,x) = u_\tau (x),
\end{equation}
where $u(t,x)$ represents the magnetization density in $x\in\mathbb{R}$ at time $t\in[\tau,+\infty)$; $\beta>0$ the inverse
temperature of the Ising system. The kernel of the convolution $J$ is a non negative, even function on $\mathbb{R}^N$ of class $\mathcal{C}^1$ with integral equal to 1 supported in the ball centered at the origin of radius 1. The kernel $J$ is related to the (long range) coupling of the spin-spin interaction; and $g$ is a globally Lipschitz continuous function of class $\mathcal{C}^1$ on $\mathbb{R}$ with $g(0)=0$.

The function $h$ is a non negative on $\mathbb{R}$ such that there exists a constant $h^*>0$, defined implicitly by (\ref{Def-h}) below, such that
\begin{equation}\label{LimH}
0\leqslant h(t,s)< h^*,\ \ (t,s)\in\mathbb{R}^2,
\end{equation}
and represents a non constant external magnetic field,  we also will assume that $h(t,\cdot)$ is a Lispschitz function with constant $\ell_h>0$, and for all $t\in\mathbb{R}$
\begin{equation}
h(t,\cdot)\in\mathcal{C}^1(\mathbb{R})\quad\mbox{with}\quad h(t,0)=0.
\end{equation}

It is interesting to note that if we take $g(t)\equiv t$, $\beta=1$ and $h(t)\equiv0$, then the linear map $Au=-u+J*u$ shares some properties with the Laplace operator, such as a form of maximum principle (see Theorem 2.1 and 2.2 in \cite{CD}). One can also see that $A$ is a nonpositive operator on $L^2(\mathbb{R}^N)$ by taking Fourier transforms since  $\widehat{J}(\xi)=\int_{\mathbb{R}^N}e^{i\xi\cdot x}J(x)dx$ is real and bounded by $1$.

Let $\mathbb{R}_\tau=[\tau,+\infty)$ for any $\tau\in\mathbb{R}$, $\rho$ a positive continuous function on $\mathbb{R}^N$ with norm equal to $1$ on $L^{1}(\mathbb{R}^N)$, and $\Omega\subset\mathbb{R}^N$ be an open set (not necessary bounded). Given $p\in[1,+\infty)$, the Banach space $L^p(\Omega,\rho)$ is defined by
\begin{equation}\label{Def_Space}
L^p(\Omega,\rho):=\Big\{u\in L_{loc}^1(\Omega);\ \int_{\Omega}\rho(x)|u(x)|^pdx<+\infty\Big\},
\end{equation}
with the norm
\begin{equation*}
\|u\|_{L^p(\Omega,\rho)}:=\Big(\int_{\Omega}\rho(x)|u(x)|^pdx\Big)^{1/p}.
\end{equation*}
We notice that the constant functions are on $L^p(\Omega,\rho)$ and $u(x)\equiv1$ has norm $1$. The corresponding higher-order weighted Sobolev space $W^{\mathfrak{m},p}(\Omega,\rho)$, $\mathfrak{m}\in\mathbb{N}$, is the space of functions $u\in L^p(\Omega,\rho)$ whose distributional derivatives up to order $\mathfrak{m}$ are also in $L^p(\Omega,\rho)$, with norm
\[
\|u\|_{W^{\mathfrak{m},p}(\Omega,\rho)}:=\Big(\sum_{|\alpha|\leqslant \mathfrak{m}}\|D^\alpha u\|^p_{L^p(\Omega,\rho)}\Big)^{1/p},\quad D^\alpha u=\partial^{|\alpha|}u/\partial x_1^{\alpha_1}\partial x_2^{\alpha_2}\cdots \partial x_N^{\alpha_N},\ \ \alpha\in\mathbb{N}^N.
\]

We observe that if the Lebesgue measure of $\Omega$ is finite and $1\leqslant p<q<+\infty$, then $ L^q(\Omega,\rho)\subset  L^p(\Omega,\rho)$ and $\|u\|_{L^p(\Omega,\rho)}\leqslant C\|u\|_{L^q(\Omega,\rho)}$ for any $u\in L^q(\Omega,\rho)$ where $C>0$ dependent of $\Omega$, $p$ and $q$.

We will see that the integral representation of the global solutions of (\ref{1.1}) in  $L^p(\Omega,\rho)$ is give by
\begin{equation}\label{Solution}
u(t,x)=e^{-(t-\tau)}u_\tau(x)+\int_\tau^te^{-(t-s)}g(\beta(J*u)(s,x) +\beta h(s,u(s,x)))ds,
\end{equation}
for all $(t,x)\in \mathbb{R}_\tau\times\mathbb{R}^N$, where $\mathbb{R}_\tau=[\tau,+\infty)$.

The outline of the paper is as follows. In section 2, we define the functional spaces and we recall some definitions of the theory of pullback attractors; in Section 3 we show the well posedness of (\ref{1.1}) on Banach spaces $L^p(\Omega,\rho)$, $p\in(1,+\infty)$ (see (\ref{Def_Space}) below), in the following sections we prove the existence of pullback attractors and study its properties. The existence of the attractor in $L^p(\mathbb{R}^N,\rho)$ for the nonlinear evolution process $S(t,\tau)u(\tau,x):=u(t,x)$, where $u(t,x)$ is given by (\ref{Solution}), is proved in the Section 4 following the ideas of \cite{Pereira} and using some estimates obtained in section 3. In Section 5, using similar arguments to \cite{carabalho} and \cite{carabalho1}, we show that the pullback attractor is a bounded set in $W^{1,p}(\mathbb{R}^N)$ and $\mathcal{C}^1(\mathbb{R}^N)$. Finally, in Section 6 we prove the upper semicontinuity of the attractors with respect to functional parameter $h(t)$ using standard techniques based on the continuity of the processes.

\section{Notations and Definitions}\label{Sec2}

For convenience of the reader, we remember the definition of \textit{nonlinear evolution process} (or non-autonomous dynamical systems) generated by problem of the type (\ref{1.1}) and pullback attractors, see \cite{CLR},\cite{Chep}, \cite{Kl}, \cite{KS} and \cite{Sell}.

An evolution process in $L^2(\Omega,\rho)$ is a family of maps $\{S(t,\tau);\ t\in\mathbb{R}_\tau,\ \tau\in\mathbb{R}\}$ from $L^2(\Omega,\rho)$ into itself with the following properties:
\begin{itemize}
\item $S(t,t)=I$, for all $t\in\mathbb{R}$,
\item $S(t,\tau)=S(t,s)S(s,\tau)$, for all $t\geqslant s\geqslant\tau$,
\item The map $\{(t,\tau)\in\mathbb{R}^2;\  t\in\mathbb{R}_\tau\}\times L^2(\Omega,\rho)\ni (t,\tau,x)\mapsto S(t,\tau)x\in L^2(\Omega,\rho)$ is continuous.
\end{itemize}

In the particular, if each $S(t,\tau)$ is linear, we say that $\{S(t,\tau);\ t\in\mathbb{R}_\tau,\ \tau\in\mathbb{R}\}$ is a linear evolution process.

A \textit{globally-defined solution} (or simply a global solution) of the nonlinear evolution process $\{S(t,\tau);\ t\in\mathbb{R}_\tau,\ \tau\in\mathbb{R}\}$ generated by equations of the type (\ref{1.1}) through $\xi_\tau\in L^2(\Omega,\rho)$ is a function $\xi:\mathbb{R}\to L^2(\Omega,\rho)$ such that $\xi(\tau)=\xi_\tau$ and for all $t\geqslant\tau$ we have $S(t,\tau)\xi(\tau)=\xi(t)$.

A subset $B$ of $L^2(\Omega,\rho)$ \textit{pullback absorbs bounded subsets of $L^2(\Omega,\rho)$ under} $\{S(t,\tau);\ t\in\mathbb{R}_\tau,\ \tau\in\mathbb{R}\}$ if $B$ pullback absorbs all bounded subsets at time $t\in\mathbb{R}$ under the process $\{S(t,\tau);\ t\in\mathbb{R}_\tau,\ \tau\in\mathbb{R}\}$, for each $t\in\mathbb{R}$, i.e., for each bounded subset $D$ of $L^2(\Omega,\rho)$, there exists $\tau_0=\tau_0(t,B)$ with $S(t,\tau)D\subset B$ for any $\tau\leqslant\tau_0$.

A family of sets $\{K(t);\ t\in\mathbb{R}\}$ \textit{pullback attracts bounded subsets of $L^2(\Omega,\rho)$ under} $\{S(t,\tau);\ t\in\mathbb{R}_\tau,\ \tau\in\mathbb{R}\}$ if $K(t)$ pullback attracts all bounded subsets at $t$ under the process $\{S(t,\tau);\ t\in\mathbb{R}_\tau,\ \tau\in\mathbb{R}\}$, for each $t\in\mathbb{R}$, i.e., for each bounded subset $C$ of $ L^2(\Omega,\rho)$
\[
\lim_{\tau\to-\infty}\mathrm{dist}(S(t,\tau)C,K(t))=0,
\]
where $\mathrm{dist}(\cdot,\cdot)$ denotes the Hausdorff semi-distance
\[
\mathrm{dist}(A,B)=\sup_{a\in A}\inf_{b\in B}|a-b|.
\]

The \textit{pullback omega-limit set} at time $t$ of a subset $B$ of $L^2(\Omega,\rho)$ is defined by
\begin{equation}\label{Def_Om}
\omega_\wp(B,t):=\bigcap_{s\leqslant t}\overline{\bigcup_{\tau\leqslant s}S(t,\tau)B}.
\end{equation}

In the sequel we introduce the concept of pullback attractor (see \cite{Kl} and \cite{KS} for more details).

A family $\{\mathcal{A}(t);\ t\in\mathbb{R}\}$ of compact subsets of $L^2(\Omega,\rho)$ is said to be the \textit{pullback attractor} for the evolution process $\{S(t,\tau);\ t\in\mathbb{R}_\tau,\ \tau\in\mathbb{R}\}$ if it is invariant, i.e., $S(t,\tau)\mathcal{A} (\tau) = \mathcal{A} (t )$ for all $t\in\mathbb{R}_\tau$, pullback attracts bounded subsets of $L^2(\Omega,\rho)$, and is minimal, that is, if there is another family of closed sets $\{C(t);\ t\in\mathbb{R}\}$ which pullback attracts bounded subsets of $L^2(\Omega,\rho)$, then $\mathcal{A}(t) \subset C(t)$, for all $t\in\mathbb{R}$. In \cite{Chep} the sets $\mathcal{A}(t)$ are referred to as kernel section.

In the non-autonomous case, the definition of pullback attractor has the same characterization
as the union of all globally-defined bounded orbits
$$
\{\mathcal{A}(t); t\in \mathbb{R}\}=\{\xi(t);\ \xi:\mathbb{R}\to L^2(\Omega,\rho)\ \hbox{is\ bounded\ and} \ S(t,\tau)\xi(\tau)=\xi(t),\ t\in\mathbb{R}_\tau,\ \tau\in\mathbb{R}\}
$$
of the autonomous case.

\section{Estimates and well-posedness in $L^2(\mathbb{R}^N,\rho)$} \label{Sec2}

In order to obtain  well posedness of (\ref{1.1}) on $L^2(\mathbb{R}^N,\rho)$, we initially  consider the following non-autonomous ODE on $L^2(\mathbb{R}^N,\rho)$
\begin{equation}\label{CP}
\partial_tu = f(t,u)\quad\mbox{and}\quad u(\tau,x)=u_{\tau}(x)
\end{equation}
where the map
\begin{equation}\label{FFF}
f(t,u)=-u+g(\beta(J*u)+\beta h(t,u))
\end{equation}
is defined on $\mathbb{R}_\tau\times L^2(\mathbb{R}^N,\rho)$.

We consider the linear operator $A$ defined by $Au:=J*u$ in $L^2(\mathbb{R}^N,\rho)$. In the analysis that follows, we use some estimates that the operator $A$ shares with the operator $A$ defined in $L^2(\mathbb{R}^N)$.

\begin{lemma}\label{Lem1}
The following statements are holds.
\item[(a)] Assume that there exists a constant $K>0$ such that $\sup\{\rho(x);\ |x-y|\leqslant1\}\leqslant K\rho(y)$ for all $y\in\mathbb{R}^N$. Then, for any $u\in L^2(\mathbb{R}^N,\rho)$
\[
\|J*u\|_{L^2(\mathbb{R}^N,\rho)}\leqslant K^{1/2}\|u\|_{L^2(\mathbb{R}^N,\rho)}
\]
and
\[
\|J'*u\|_{L^2(\mathbb{R}^N,\rho)}\leqslant  K^{1/2}\|u\|_{L^2(\mathbb{R}^N,\rho)}.
\]
\item[(b)] There exists $C>0$ such that for any $u,v\in L^2(\mathbb{R}^N,\rho)$
\[
|J*(u+v)(x)|\leqslant C\|u+v\|_{L^2(\mathbb{R}^N,\rho)}
\]
for all $x\in\mathbb{R}^N$.
\end{lemma}
$\proof$ (a) By H\"{o}lder's inequality, we have
\begin{eqnarray*}
|(J*u)(x)|&\leqslant& \int_{\mathbb{R}^N}|J(x-y)||u(y)|dy\\
&=& \int_{\mathbb{R}^N}|J(x-y)|^{1/2}|J(x-y)|^{1/2}|u(y)|dy\\
&\leqslant& \Big(\int_{\mathbb{R}^N}|J(x-y)|dy\Big)^{(p-1)/p}\Big(\int_{\mathbb{R}^N}|J(x-y)||u(y)|^pdy\Big)^{1/2}\\
&\leqslant& \Big(\int_{\mathbb{R}^N}|J(x-y)||u(y)|^pdy\Big)^{1/p}.
\end{eqnarray*}
By Fubini's theorem
\begin{eqnarray*}
\|J*u\|_{L^2(\mathbb{R}^N,\rho)}^p&=&\int_{\mathbb{R}^N}\rho(x)|(J*u)(x)|^p dx\\
&\leqslant& \int_{\mathbb{R}^N}\rho(x)\Big(\int_{\mathbb{R}^N}|J(x-y)||u(y)|^pdy \Big)dx\\
&=& \int_{\mathbb{R}^N}\Big(\int_{\mathbb{R}^N}\rho(x)|J(x-y)|dx \Big)|u(y)|^pdy
\end{eqnarray*}
and using the characteristic function, $\chi_{B[y;1]}$, of the ball $B[y;1]$ in $\mathbb{R}^N$,  we get
\begin{eqnarray*}
\|J*u\|_{L^2(\mathbb{R}^N,\rho)}^p&=&\int_{\mathbb{R}^N}\Big(\int_{\mathbb{R}^N}\rho(x)|J(x-y)|\chi_{B[y;1]}(x)dx \Big)|u(y)|^pdy\\
&\leqslant& K\int_{\mathbb{R}^N}\rho(y)|u(y)|^pdy\\
&=& K\|u\|_{L^2(\mathbb{R}^N,\rho)}^p
\end{eqnarray*}
which  concludes the item (a).

The second estimate of (a) follows from same arguments.

Now, we will show (b), using the characteristic function, $\chi_{B[0;1]}$, of the ball $B[0;1]$ in $\mathbb{R}^N$, we notice that for each $x\in\mathbb{R}^N$
\begin{eqnarray*}
|J*(u+v)(x)|&\leqslant&\int_{\mathbb{R}^N}|J(x-y)||(u+v)(y)|dy\\
&\leqslant&\int_{\mathbb{R}^N}|J(\xi)||(u+v)(x-\xi)|d\xi\\
&\leqslant&\int_{\mathbb{R}^N}\chi_{B[0;1]}(\xi)|J(\xi)||(u+v)(x-\xi)|d\xi\\
&\leqslant&\|J\|_{L^\infty(\mathbb{R}^N)}\rho_1^{-1}\int_{\mathbb{R}^N}\chi_{B[0;1]}(\xi)\rho(\xi)|(u+v)(x-\xi)|d\xi\\
&=&\|J\|_{L^\infty(\mathbb{R}^N)}\rho_1^{-1}\int_{\mathbb{R}^N}\chi_{B[x;1]}(\eta)\rho(\eta)|(u+v)(\eta)|d\eta
\end{eqnarray*}
where $\rho_1=\inf\{\rho(y); |y|\leqslant 1\}>0$.

Finally, since $L^p(B[z;1],\rho)\subset L^1(B[z;1],\rho)$ and $\|\omega\|_{L^1(B[z;1],\rho)}\leqslant\|\omega\|_{L^p(B[z;1],\rho)}$, for any $z\in\mathbb{R}^N$ and $\omega\in L^p(B[z;1],\rho)$, we get
\begin{eqnarray*}
|J*(u+v)(x)|&\leqslant&\|J\|_{L^\infty(\mathbb{R}^N)}\rho_1^{-1}\Big(\int_{\mathbb{R}^N}\chi_{B[x;1]}(\eta)\rho(\eta)|(u+v)(\eta)|^pd\eta\Big)^{1/p}\\
&\leqslant&\|J\|_{L^\infty(\mathbb{R}^N)}\rho_1^{-1}\Big(\int_{\mathbb{R}^N}\rho(\eta)|(u+v)(\eta)|^pd\eta\Big)^{1/p}\\
&=&\|J\|_{L^\infty(\mathbb{R}^N)}\rho_1^{-1}\|u+v\|_{L^2(\mathbb{R}^N,\rho)}
\end{eqnarray*}
so the proof is complete. $\blacksquare$

\begin{remark}[\cite{Pereira}]
The hypothesis $\sup\{\rho(x);\ |x-y|\leqslant1\}\leqslant K\rho(y)$ of the Lemma \ref{Lem1} in (a) is verified, for instance, if $N=1$ and $\rho(x)=\pi^{-1}(1+x^2)^{-1}$, with $K=3$.
\end{remark}

We will show that $f$ is a globally Lipschitz continuous function on $L^2(\mathbb{R}^N,\rho)$ with respect to the second variable.

\begin{proposition}\label{Prop well posed}
Assume that $g$ is globally Lipschitz continuous in $\mathbb{R}$ with constant $\ell_g>0$. For each $t\in\mathbb{R}$, the map $f(t,\cdot)$ given in (\ref{FFF}) is globally Lipschitz continuous in $L^2(\mathbb{R}^N,\rho)$ with
\[
\|f(t,u)-f(t,v)\|_{L^2(\mathbb{R}^N,\rho)}\leqslant (1+\ell_g\beta K^{1/p}+\beta\ell_h)\|u-v\|_{L^2(\mathbb{R}^N,\rho)}
\]
where $K>0$ is such that $\sup\{\rho(x);\ |x-y|\leqslant1\}\leqslant K\rho(y)$ for all $y\in\mathbb{R}^N$ (as well as in Lemma \ref{Lem1}).
\end{proposition}
\noindent$\proof$ Since $J$ is bounded and compact supported, $(J*u)(x)$ is well defined for $u\in L^1_{loc}(\mathbb{R}^N)$. Since $g$ is globally bounded by constant $a>0$ it follows that $f(t,u)\in L^2(\mathbb{R}^N,\rho)$ if $(t,u)\in\mathbb{R}\times L^2(\mathbb{R}^N,\rho)$. Now, from Lemma \ref{Lem1}, we have
\begin{eqnarray*}
&&\|f(t,u)-f(t,v)\|_{L^2(\mathbb{R}^N,\rho)}\\
&\leqslant&\|u-v\|_{L^2(\mathbb{R}^N,\rho)}+\| g(\beta(J*u)+\beta h(t,u))-g(\beta(J*v) +\beta h(t,v))\|_{L^2(\mathbb{R}^N,\rho)}\\
&\leqslant &\|u-v\|_{L^2(\mathbb{R}^N,\rho)}+\ell_g\beta\|J*(u-v)\|_{L^2(\mathbb{R}^N,\rho)}+\beta\ell_h\|u-v\|_{L^2(\mathbb{R}^N,\rho)}\\
&\leqslant &(1+\ell_g\beta K^{1/p}+\beta\ell_h)\|u-v\|_{L^2(\mathbb{R}^N,\rho)}
\end{eqnarray*}
as claimed. $\blacksquare$

From Proposition \ref{Prop well posed} and basic theory of ODE's in Banach spaces it follows that, for any $u_\tau\in L^2(\mathbb{R}^N,\rho)$, the Cauchy problem (\ref{CP}) has a unique local solution in\linebreak $\mathcal{C}([\tau,s(u_\tau)], L^2(\mathbb{R}^N,\rho))\cap \mathcal{C}^1((\tau,s(u_\tau)], L^2(\mathbb{R}^N,\rho))$ for some $s(u_\tau)>0$  which is continuous with respect to $u_\tau$. By standard arguments, using the variation of constants formula and Gronwall's inequality, it follows that these solutions are actually globally defined, i.e., $s(u_\tau)=\infty$ for any $u_\tau$.

The natural notation for the global solution of the Cauchy problem (\ref{CP}) is $u(t,\tau, x;u_\tau)$. In this paper for simplicity of notation, we use $u(t,x)$ to denote the global solution.

\section{Existence of pullback attractors} \label{Sec3}

In this section, we will prove that $S(t,\tau)u(\tau,x):=u(t,x)$ ($t\in\mathbb{R}_\tau$, $\tau\in\mathbb{R}$), where $u(t,x)$ denote the global solution of the Cauchy problem (\ref{CP}), provides an infinite-dimensional non-autonomous dynamical system in $L^2(\mathbb{R}^N,\rho)$ that has a pullback attractor $\{\mathcal{A}(t);\ t\in\mathbb{R}\}$.

The next result is a extension of Lemma 3 of \cite{Pereira}.

\begin{lemma}\label{Lema 3.1}
Assume the same hypotheses of Proposition \ref{Prop well posed}. If $g$ is globally bounded by a constant $a>0$, then the ball $B(0;a+\epsilon)$ is a pullback absorbing for the evolution process $S(t,\tau)$ generated by (\ref{CP}) in $L^2(\mathbb{R}^N,\rho)$, for any $\epsilon>0$.
\end{lemma}
\noindent $\proof$ Let $u(t,x)$ be the solution of (\ref{CP}) with initial condition $u(\tau,x)\in B$, where $B$ is a bounded subset of $L^2(\mathbb{R}^N,\rho)$, namely
\begin{equation}\label{aaa}
u(t,x)=e^{-(t-\tau)}u(\tau,x)+\int_{\tau}^{t}e^{-(t-s)}g(\beta(J*u)(s,x)+\beta h(s,u(s,x)))ds.
\end{equation}
We observe that
\begin{eqnarray*}
&&\|\int_{\tau}^t e^{-(t-s)} g(\beta(J*u)(s,x)+\beta h(s,u(s,x)))ds\|^p_{L^2(\mathbb{R}^N,\rho)}\\
&\leqslant&\int_{\mathbb{R}^N}\rho(x)\left(\int_{\tau}^t e^{-(t-s)} |g(\beta(J*u)(s,x)+\beta h(s,u(s,x)))|ds\right)^p dx.
\end{eqnarray*}
Using the boundedness of $g$  and the fact that the norm of $\rho$ is equal to 1, it follows that
\begin{eqnarray*}
\|\int_{\tau}^t e^{-(t-s)} g(\beta(J*u)(s,x)+\beta h(s,u(s,x)))ds\|^p_{L^2(\mathbb{R}^N,\rho)}&\leqslant& a^p\int_{\mathbb{R}^N}\rho(x)\left(\int_{\tau}^t e^{-(t-s)}ds\right)^p dx\\
&\leqslant& a^p.
\end{eqnarray*}
Hence
\begin{equation}\label{aaaa}
\|\int_{\tau}^t e^{-(t-s)} g(\beta(J*u)(s,x)+\beta h(s,u(s,x)))ds\|_{L^2(\mathbb{R}^N,\rho)}\leqslant a.
\end{equation}

Now, we notice that from (\ref{aaa}), (\ref{aaaa}) and Minkoviski's inequality
\begin{equation*}
\|u(t,\cdot)\|_{L^2(\mathbb{R}^N,\rho)}\leqslant e^{-(t-\tau)}\|u_\tau\|_{L^2(\mathbb{R}^N,\rho)}+a.
\end{equation*}
Therefore,  if $\|u_\tau\|_{L^2(\mathbb{R}^N,\rho)}\neq0$ then $u(t,x)\in B(0;a+\epsilon)$ for all $\tau\leqslant \tau_0(t,B)$, where $\tau_0(t,B)=\displaystyle{\ln\left(\frac{\epsilon}{\|u_\tau\|_{L^2(\mathbb{R}^N,\rho)}}\right)+t}$.

On other hand, if $\|u_\tau\|_{L^2(\mathbb{R}^N,\rho)}=0$ then $\|u(t,\cdot)\|_{L^2(\mathbb{R}^N,\rho)}\leqslant a$, i.e., $u(t,x)\in B(0;a)$ for all $\tau\in\mathbb{R}$.

Therefore, for each $t\in\mathbb{R}$ there exists $\tau_0=\tau(t,B)\in\mathbb{R}$ such that for any $\tau\leqslant\tau_0$
\[
S(t,\tau)B\subset B(0;a+\epsilon).
\]
Then, the result follows.
$\blacksquare$

\bigskip

Our next goal is to prove that the pullback attractor is the family of pullback omega-limit \{$\omega_\wp(B(0;a+\epsilon),t);\ t\in\mathbb{R}\}$ (see (\ref{Def_Om})). Next, we established a result that is an extension of Lemma 4 of \cite{Pereira}. It will be used to prove the compactness of the sets $\omega_\wp(B(0;a+\epsilon),t)$.

\begin{lemma}\label{Lema 4}
Assume the same hypotheses of Lemma \ref{Lema 3.1}. If $g'$ is globally Lipschitz continuous with constant $k_{1}>0$ and $k_2=|g'(0)|>0$, then for any $\eta>0$ and $t\in\mathbb{R}$ there exists $\tau_\eta\leqslant t$ such that $S(t,\tau_\eta)B(0;a+\epsilon)$ has a finite covering by balls of $L^2(\mathbb{R}^N,\rho)$ with radius smaller than $\eta$.
\end{lemma}
\noindent $\proof$ From Lemma \ref{Lema 3.1} it follows that $S(t,\tau)B(0;a+\epsilon)\subset B(0;a+\epsilon)$, for any $t\in\mathbb{R}_\tau$, $\tau\in\mathbb{R}$. Given $u_\tau\in B(0;a+\epsilon)$, we consider the non-autonomous coupled system
\begin{equation}\label{CoSist1}
\left\{
\begin{array}{llll}
\partial_t v(t,x)=-v(t,x)\\
\partial_t w(t,x)=-w(t,x)+g(\beta(J*(v+w)(t,x))+\beta h(t,(v+w)(t,x)))
\end{array}
\right.
\end{equation}
with initial conditions at $t=\tau$, $\tau\in\mathbb{R}$
\begin{equation}\label{CoSist2}
\left\{
\begin{array}{ll}
\displaystyle v(\tau,x)=u_\tau(x),\\
w(\tau,x)=0
\end{array}
\right.
\end{equation}
which is a way to rewrite the ODE (\ref{CP}) on $L^p(\mathbb{R}^N, \rho)$.

Note that if $(v, w)$ is the solution of (\ref{CoSist1})-(\ref{CoSist2}) in $L^2(\mathbb{R}^N,\rho)\times L^2(\mathbb{R}^N,\rho)$ then $u = v + w$ is a solution of (\ref{CP}) in $L^2(\mathbb{R}^N,\rho)$ with $u(\tau) = u_\tau$. Conversely, any solution $u$ of (\ref{CP})-(\ref{CP_1}) in $L^2(\mathbb{R}^N,\rho)$ can be written as $u = v + w$, with $(v, w)$ the solution of (\ref{CoSist1})-(\ref{CoSist2}) in $L^2(\mathbb{R}^N,\rho)\times L^2(\mathbb{R}^N,\rho)$.

Since $v(t,x) = e^{-(t-\tau)}u_\tau(x)$, given $\eta>0$, we may find $\tau_\eta\in\mathbb{R}$ such that if $t\geqslant\tau_\eta\geqslant\tau$ then $\|v(t,\cdot)\|_{L^2(\mathbb{R}^N,\rho)}\leqslant\eta/2$, for any $u_\tau\in B(0;a+\epsilon)$.

By variation of constants formula
\begin{equation*}\label{EDR}
w(t,x)=\int_{\tau}^{t}e^{-(t-s)}g(\beta(J*u)(s,x)+\beta h(s,u(s,x)))ds,\quad u = v + w,
\end{equation*}
and therefore, by boundedness of $g$, we obtain
\begin{equation}\label{9InRef}
|w(t,x)|\leqslant a,\ \forall\ (t,x)\in \mathbb{R}_\tau\times \mathbb{R}^N.
\end{equation}

From now on, we write the solution $w$ of the following way
\begin{equation*}
w(t,\cdot)=w(t,\cdot)\chi_{B(0;R)}+w(t,\cdot)(1-\chi_{B(0;R)})
\end{equation*}
where $\chi_{B(0;R)}$ denotes the characteristic function of the ball $B(0;R)$ and $R>0$ is a constant to be chosen.

By item (a) of the Lema \ref{Lem1}, if $u\in B(0;a+\epsilon)$, $R>0$ and $x\in B(0;R)$, then
\begin{eqnarray*}
|(J*u)(s,x)|^p&\leqslant&\int_{\mathbb{R}^N}|J(x-y)||u(s,y)|^pdy\\
&=&\int_{\mathbb{R}^N}|J(y)||u(s,y)|^p\chi_{B(x;R+1)}(y)dy.
\end{eqnarray*}

Let $\rho_{R+1}>0$ be the infimum of the set $\{\rho(y);\ |y-x|\leqslant R+1\}$. Thus
\begin{eqnarray*}
|(J*u)(s,x)|^p&\leqslant&\frac{1}{\rho_{R+1}}\int_{\mathbb{R}^N}\rho(y)|J(y)||u(s,y)|^p\chi_{B(x;R+1)}(y)dy\\
&\leqslant&\frac{1}{\rho_{R+1}}\|J\|_{L^\infty(\mathbb{R}^N)}\int_{\mathbb{R}^N}\rho(y)|u(s,y)|^p\chi_{B(x;R+1)}(y)dy\\
&\leqslant&\frac{1}{\rho_{R+1}}\|J\|_{L^\infty(\mathbb{R}^N)}\int_{\mathbb{R}^N}\rho(y)|u(s,y)|^pdy\\
&\leqslant&\frac{1}{\rho_{R+1}}\|J\|_{L^\infty(\mathbb{R}^N)}\|u(s,\cdot)\|_{L^2(\mathbb{R}^N,\rho)}^p.
\end{eqnarray*}

Since $S(t,\tau)B(0;a+\epsilon)\subset B(0;a+\epsilon)$, for any $t\geqslant\mathbb{R}_\tau$, $\tau\in\mathbb{R}$, we have
\begin{equation}\label{Cita2}
|(J*u)(x,s)|^p\leqslant C_1(\varepsilon)
\end{equation}
where $C_1(\varepsilon)=\frac{2^{p}\|J\|_{L^\infty(\mathbb{R}^N)}}{\rho_{R+1}}(a^{p}+\varepsilon^p)>0$.

Using similar arguments we get
\begin{equation}\label{Cita3}
|( J'*u)(x,s)|^p\leqslant C_2(\varepsilon),
\end{equation}
where $C_2(\varepsilon)=\frac{2^{p}\|J'\|_{L^\infty(\mathbb{R}^N)}}{\rho_{R+1}}(a^{p}+\varepsilon^p)>0$.

Since $g'$ is  globally Lipschitz continuous on $\mathbb{R}$ with constant $k_1>0$, using (\ref{LimH}) we find
\begin{eqnarray*}
&&\left|\partial_x w(t,x)\right|\\
& \leqslant & \beta
\int_{\tau}^{t}e^{-(t-s)}|g'(\beta (J*u)(s,x)+\beta h(s,u(s,x)))||(J'*u)(s,x)+\partial_2 h(t,u(s,x))\partial_xu(s,x) |ds\nonumber\\
&\leqslant& \beta^2k_1
\int_{\tau}^{t}e^{-(t-s)}|(J*u)(x,s)+\beta h(s,u(s,x))||( J'*u)(s,x)+\partial_2 h(t,u(s,x))\partial_xu(s,x)|ds\nonumber\\
&+&\beta k_2h^*
\int_{\tau}^{t}e^{-(t-s)}|( J'*u)(s,x)+\partial_2 h(t,u(s,x))\partial_xu(s,x)|ds\nonumber\\
&\leqslant& \beta^2k_1
\int_{\tau}^{t}e^{-(t-s)}|(J*u)(s,x)||(J'*u)(s,x)+\partial_2 h(t,u(s,x))\partial_xu(s,x)|ds\nonumber\\
&+&\beta h^*(1+k_2)\int_{\tau}^{t}e^{-(t-s)}|(J'*u)(s,x)+\partial_2 h(t,u(s,x))\partial_xu(s,x)|ds\nonumber
\end{eqnarray*}
for any $t \in\mathbb{R}_\tau$, $\tau\in\mathbb{R}$, $x\in B(0;R)$.

Hence
\begin{eqnarray*}
\left|\partial_x w(t,x)\right| &\leqslant& \beta^2k_1C_1(\varepsilon)C_2(\varepsilon)+\beta h^*(1+k_2)C_2(\varepsilon)   \nonumber\\
&+&\beta^2k_1h_1\int_{\tau}^{t}e^{-(t-s)} |(J*u)(s,x)||\partial_xu(s,x)|ds\\
&+& \beta h^*h_1(1+k_2)\int_{\tau}^{t}e^{-(t-s)}|\partial_xu(s,x)|ds\nonumber
\end{eqnarray*}
for any $t \in\mathbb{R}_\tau$, $\tau\in\mathbb{R}$, $x\in B(0;R)$.

Therefore, by (\ref{Cita2}) and (\ref{Cita3}) we obtain
\begin{equation*}
\left|\partial_x w(t,x)\right| \leqslant  \beta^{2}k_1C_1(\varepsilon)C_2(\varepsilon)+\beta( k_2+h^*)C_2(\varepsilon).
\end{equation*}

Let $R>0$ be chosen such that
\[
\int_{\mathbb{R}^N}\rho(x)(1-\chi_{B(0;R)}(x))dx\leqslant\frac{\eta^p}{4^pa^p}.
\]
Then, by (\ref{9InRef}) we get
$$
\|w(t,\cdot)(1-\chi_{B(0;R)})\|^p_{L^2(\mathbb{R}^N,\rho)}\leqslant\int_{\mathbb{R}^N}\rho(x)(1-\chi_{B(0;R)}(x))|w(t,x)|^pdx\leqslant\frac{\eta}{4}.
$$

Moreover, by (\ref{9InRef}) the function $w(t,\cdot)\chi_{B(0;R)}$ is bounded in $W^{1,p}(B(0;R),\rho)$ (by a constant independent of $u\in B(0;a+\epsilon)$) and, therefore the set $\{w(t,x);\ x\in B(0;R)\}$ with $w(\tau,\cdot)\in B(0;a+\epsilon)$ is a compact subset of $L^p(B(0;R),\rho)$ for any $t\in\mathbb{R}_\tau$ and, thus, it can be covered by a finite number of balls with radius smaller than $\eta/4$.

Therefore, since $u(t,\cdot) $ is the solution of the system (\ref{CoSist1})-(\ref{CoSist2}) in $L^2(\mathbb{R}^N,\rho)$ we can be writer as
\[
u(t,\cdot) = v(t,\cdot) + w(t,\cdot)\chi_{B(0;R)}+w(t,\cdot)(1-\chi_{B(0;R)})
\]
it follows that $S(t,\tau_\eta)B(0;a+\epsilon)$ has a finite covering by balls of $L^2(\mathbb{R}^N,\rho)$ with radius smaller than $\eta$. $\blacksquare$

\begin{theorem}\label{Theorem 3.2}
Assume the same hypotheses of Lemma \ref{Lema 4}. The family of sets $\mathcal{A}(t)=\omega_\wp(B(0;a+\epsilon),t)$ is a pullback attractor for the process $S(t,\tau)$ generated by (\ref{CP}) in $L^p(\mathbb{R}^N)$. Moreover, $\bigcup_{t\in\mathbb{R}}\mathcal{A}(t)$ is contained in the ball of radius $a>0$.
\end{theorem}
\noindent $\proof$ From Lemma \ref{Lema 3.1}, it follows that $\mathcal{A}(t)$ is contained in the ball of radius $a$ in $L^p(\mathbb{R}^N)$ for any $t\in\mathbb{R}$. Also, since $\mathcal{A}(t)$ is positively invariant by the process, i.e., $S(t,\tau)\mathcal{A}(\tau)=\mathcal{A}(t)$ for all $t\geqslant\tau$, $\tau\in\mathbb{R}$, it follows that $\mathcal{A}(t)\subset S(t,\tau)B(0;a+\epsilon)$ for any $t\in\mathbb{R}_\tau$, and then, from Lemma \ref{Lema 4}, we obtain that the measure of noncompactness of $\mathcal{A}(t)$ is zero, for any $t\in\mathbb{R}$. Thus, $\mathcal{A}(t)$ is relatively compact and, being closed, also compact.

Finally, it remains to prove that $\mathcal{A}(t)$ pullback attracts bounded subset of $L^2(\mathbb{R}^N,\rho)$ at time $t$. If $D$ is a bounded subset of $L^2(\mathbb{R}^N,\rho)$ then $S(t,\bar{\tau})D\subset B(0;a+\epsilon)$ for $\bar{\tau}$ rather small and, therefore, $\omega_\wp(D,t)\subset\omega_\wp(B(0;a+\epsilon),t)=\mathcal{A}(t)$ for any $t\in\mathbb{R}$. $\blacksquare$

\section{Regularity of the attractors}\label{Sec5}

In this section, we will show that the pullback attractor is contained a fixed bounded subset of the Banach spaces $W^{1,p}(\mathbb{R}^N,\rho)$ and $\mathcal{C}^1(\mathbb{R}^N)$.

First, we will prove that the attractor is a bounded subset of $W^{1,p}(\mathbb{R}^N,\rho)$. Since the attractor can be written as the set of all global bounded solutions, if $u(t,x)$ is a solution of (\ref{CP}) in $\mathcal{A}(t)$ for all $t\in\mathbb{R}$, then we obtain, letting $\tau\to-\infty$
\begin{equation}\label{Iddd}
u(t,x)=\int_{-\infty}^{t}e^{-(t-s)}g(\beta(J*u)(s,x)+\beta h(s,u(s,x)))ds.
\end{equation}
Due to Lemma \ref{Lema 4} $u(t,x)\in C(\mathbb{R}\times\mathbb{R}^N, W^{1,p}(\mathbb{R}^N,\rho))$ and $\partial_t u(t,x)\in C(\mathbb{R}\times\mathbb{R}^N, L^2(\mathbb{R}^N,\rho))$. Therefore $\mathcal{A}(t)=S(t,\tau)\mathcal{A}(\tau)\in W^{1,p}(\mathbb{R}^N,\rho)$ for all $t\in\mathbb{R}$.

Furthermore, we have

\begin{theorem}\label{Theorem6}
For each $t\in\mathbb{R}$, the set $\mathcal{A}(t)$ is bounded in $\mathcal{C}^1$.
\end{theorem}
\noindent $\proof$ If $u(t,x)$ is the solution of (\ref{CP}) in $\mathcal{A}(t)$, from (\ref{Iddd}) we have
\[
u(t,x)=\int_{-\infty}^{t}e^{-(t-s)}g(\beta(J*u)(s,x)+\beta h(s,u(s,x)))ds.
\]
The equality above is in the sense of $L^2(\mathbb{R}^N,\rho)$ but, since the right-hand side is regular as $J$ we have
\begin{equation}\label{Desig}
|u(t,x)|\leqslant a\int_{-\infty}^{t}e^{-(t-s)}\leqslant a,\ \forall (t,x)\in \mathbb{R}_\tau\times \mathbb{R}^N.
\end{equation}
From (\ref{Desig}) we obtain
\begin{equation}\label{Desigualdade}
|(J'* u)(t,x) |\leqslant a\|J'\|_{L^1(\mathbb{R}^N)}\ \mathrm{and}\  |(J'* u)(t,x) |\leqslant a\|J'\|_{L^1(\mathbb{R}^N)},\ \forall (t,x)\in \mathbb{R}_\tau\times \mathbb{R}^N.
\end{equation}

Differential in (\ref{Iddd}) with respect to $x$, we obtain for $t\geqslant\tau$
\[
\partial_xu(t,x)=\beta\int_{-\infty}^te^{-(t-s)}g'(\beta (J*u)(s,x)+\beta h(s,u(s,x)))((J'*u)(s,x)+\partial_xh(s,u(s,x)))ds
\]
which is well defined by arguments entirely similar to the ones used in the proof of Lemma \ref{Lema 4}.

Since $g'$ is  globally Lipschitz continuous on $\mathbb{R}$ with constant $k_1>0$ we find
\begin{eqnarray*}
&&|\partial_x u(t,x)|\\
&\leqslant &\beta \int_{-\infty}^te^{-(t-s)}|g'(\beta (J*u)(s,x)+\beta h(s,u(s,x)))||(J'*u)(s,x)+\partial_2h(s,u(s,w))\partial_x u(s,x)|ds\\
&\leqslant& \beta^2k_1\int_{-\infty}^{t}e^{-(t-s)}|(J*u)(s,x)||(J'*u)(s,x)+\partial_2h(s,u(s,w))\partial_x u(s,x)|ds\\
&+&\beta^2h^*(1+k_2)\int_{-\infty}^{t}e^{-(t-s)}|(J'*u)(s,x)+\partial_2h(s,u(s,w))\partial_x u(s,x)|ds
\end{eqnarray*}
and from  (\ref{Desigualdade}) we obtain
\begin{eqnarray*}
\|\partial_x u(t,\cdot)\|_{L^2(\mathbb{R}^N,\rho)}&\leqslant &\beta^2k_1\int_{-\infty}^{t}e^{-(t-s)} a^2 \|J'\|_{L^1(\mathbb{R}^N)}\|J\|_{L^1(\mathbb{R}^N)}  ds\\
&+&\beta^2h^*(1+k_2)\int_{-\infty}^{t}e^{-(t-s)}a\|J'\|_{L^1(\mathbb{R}^N)} ds\\
&\leqslant& a\beta^2\|J'\|_{L^1(\mathbb{R}^N)}[ak_1 \|J\|_{L^1(\mathbb{R}^N)} +k_2+h^*]
\end{eqnarray*}
concluding the proof. $\blacksquare$

\section{Upper semicontinuity of the attractors}\label{Sec6}

We suppose that there exist functions $h_\epsilon:\mathbb{R}\to\mathbb{R}$ satisfying (\ref{LimH}), for any $\epsilon\in[0,1]$, and
we assume the convergence $h_\epsilon(t)\to h_0(t)$, as $\epsilon\to0^+$, uniformly on $\mathbb{R}$. From now on, we will denote as $\{S_\epsilon(t,\tau);\ t\in\mathbb{R}_\tau,\ \tau\in\mathbb{R}\}$ the process associated with the problem (\ref{CP})-(\ref{CP_1}) with $h=h_\epsilon$. We prove the upper semicontinuity of the pullback attractors for (\ref{1.1})-(\ref{1.1-IC}) as $\epsilon\to0^+$, i.e., we show that
\[
\lim_{\epsilon\to0^+}\mathrm{dist}(\mathcal{A}_\epsilon(t),\mathcal{A}_0(t))=0
\]
where $\{A_\epsilon(t);\ t\in\mathbb{R}\}$ denotes the pullback attractor of $S_\epsilon(t,\tau)$ on $L^2(\mathbb{R}^N,\rho)$, for any $\epsilon\in[0,1]$.

\begin{theorem}\label{The01}
 Let $\{S_\epsilon(t,\tau);\ t\in\mathbb{R}_\tau,\ \tau\in\mathbb{R}\}$ as above. For each $u_0\in L^2(\mathbb{R}^N,\rho)$, we have
$$
\|S_\epsilon(t,\tau)u_0-S_0(t,\tau)u_0\|_{L^2(\mathbb{R}^N,\rho)}\leqslant   M_1\| h_\epsilon- h_0\|_{L^\infty(\mathbb{R}^2)}e^{M_1\|J\|_{L^\infty(\mathbb{R}^N)}\rho_1^{-1}(t-\tau)}
$$
where $M_1=2^{(p+1)/p} \ell_g\beta>0$ and $\rho_1=\inf\{\rho(\xi);\ |\xi|\leqslant 1 \}>0$.
\end{theorem}
\noindent  $\proof$ Let $u_0\in L^2(\mathbb{R}^N,\rho)$ and $u^\epsilon=S_\epsilon(t,\tau)u_0$, for any $\epsilon\in[0,1]$. Then for any $(t,x)\in\mathbb{R}_\tau\times\mathbb{R}^N$
\[
(u^\epsilon -u^0)(t,x)=\int_\tau^t e^{-(t-s)}[G(s,u^\epsilon(s,x))-G(s,u^0(s,x))]ds
\]
where $G(s,u^\epsilon(s,x)):=g(\beta (J*u^\epsilon)(s,x)+\beta h_\epsilon(s,u^\epsilon(s,x)))$.

Firstly, we notice that if $g$ is  globally Lipschitz continuous on $\mathbb{R}$ with constant $\ell_g>0$, then
\begin{eqnarray}\label{Des_F}
&& \|(u^\epsilon -u^0)(t,\cdot)\|_{L^2(\mathbb{R}^N,\rho)}^p\\
 &=&\int_{\mathbb{R}^N}\rho(x)\Big|\int_\tau^t e^{-(t-s)}[G(s,u^\epsilon(s,x))-G(s,u^0(s,x))]ds\Big|^pdx\nonumber\\
 &\leqslant&\int_{\mathbb{R}^N}\rho(x)\Big(\ell_g\beta\int_\tau^t e^{-(t-s)}|J*(u^\epsilon-u^0)(s,x)| ds\nonumber\\
 & +& \ell_g\beta\int_\tau^t e^{-(t-s)}| h_\epsilon(s,u^\epsilon(s,x))- h_0(s,u^0(s,x))|ds\Big)^pdx\nonumber\\
 &\leqslant&\int_{\mathbb{R}^N}\rho(x)\Big(\ell_g\beta\int_\tau^t e^{-(t-s)}|J*(u^\epsilon-u^0)(s,x)| ds\nonumber\\
 & +& \ell_g\beta\int_\tau^t e^{-(t-s)}| (h_\epsilon- h_0)(s,u^\epsilon(s,x))|ds\nonumber\\
 &+&\ell_g\beta\int_\tau^t e^{-(t-s)}| h_0((s,u^\epsilon)-(s,u^0))(s,x)|ds\Big)^pdx\nonumber\\
 &\leqslant&\int_{\mathbb{R}^N}\rho(x)\Big(\ell_g\beta\int_\tau^t e^{-(t-s)}|J*(u^\epsilon-u^0)(s,x)| ds+\ell_g\beta\| h_\epsilon- h_0\|_{L^\infty(\mathbb{R})}\Big)^pdx.
\end{eqnarray}
Since
\[
|a_1+a_2|^{q}\leqslant 2^{q}\max\{|a_1|,|a_2|\}^{q}\leqslant 2^{q}(|a_1|^{q}+|a_2|^{q})
\]
for any $a_1,a_2\in\mathbb{R}$ and $q\in(0,+\infty)$ it follows that
\begin{eqnarray*}
 && \|(u^\epsilon -u^0)(t,\cdot)\|_{L^2(\mathbb{R}^N,\rho)}^p\nonumber\\
 &\leqslant &M\| h_\epsilon- h_0\|^p_{L^\infty(\mathbb{R})}+ M \int_{\mathbb{R}^N}\rho(x)\Big(\int_\tau^t e^{-(t-s)}|J*(u^\epsilon-u^0)(s,x)| ds\Big)^pdx
\end{eqnarray*}
where $M=2^{p} \ell_g^p\beta^p>0$.

Now, using the item (c) of Lemma \ref{Lem1}, we find
\begin{eqnarray}\label{Des1}
 && \|(u^\epsilon -u^0)(t,\cdot)\|_{L^2(\mathbb{R}^N,\rho)}^p\nonumber\\
 &\leqslant &M\| h_\epsilon- h_0\|^p_{L^\infty(\mathbb{R})}+ M\|J\|_{L^\infty(\mathbb{R}^N)}^p\rho_1^{-p} \int_{\mathbb{R}^N}\rho(x)\Big(\int_\tau^t e^{-(t-s)}\|(u^\epsilon-u^0)(s,\cdot)\|_{L^2(\mathbb{R}^N,\rho)} ds\Big)^pdx\nonumber\\
 &= &M\| h_\epsilon- h_0\|^p_{L^\infty(\mathbb{R})}+ M\|J\|_{L^\infty(\mathbb{R}^N)}^p\rho_1^{-p} \Big(\int_\tau^t e^{-(t-s)}\|(u^\epsilon-u^0)(s,\cdot)\|_{L^2(\mathbb{R}^N,\rho)} ds\Big)^p\int_{\mathbb{R}^N}\rho(x)dx\nonumber\\
 &\leqslant &M\| h_\epsilon- h_0\|^p_{L^\infty(\mathbb{R})}+ M\|J\|_{L^\infty(\mathbb{R}^N)}^p\rho_1^{-p} \Big(\int_\tau^t e^{-(t-s)}\|(u^\epsilon-u^0)(s,\cdot)\|_{L^2(\mathbb{R}^N,\rho)} ds\Big)^p\nonumber
\end{eqnarray}
where $M=2^{p} \ell_g^p\beta^p$ and $\rho_1=\inf\{\rho(\xi);\ |\xi|\leqslant 1 \}>0$.

Hence
\begin{eqnarray}\label{Des_F}
 &&\|(u^\epsilon -u^0)(t,\cdot)\|_{L^2(\mathbb{R}^N,\rho)}\nonumber\\
 &\leqslant & M_1\| h_\epsilon- h_0\|_{L^\infty(\mathbb{R})}+M_1\|J\|_{L^\infty(\mathbb{R}^N)}\rho_1^{-1}\int_\tau^t e^{-(t-s)}\|(u^\epsilon-u^0)(s,\cdot)\|_{L^2(\mathbb{R}^N,\rho)} ds.
\end{eqnarray}
where $M_1=2^{1/p}M^{1/p}=2^{(p+1)/p} \ell_g\beta>0$.

Finally, setting $\varphi(t)=e^t\|(u^\epsilon -u^0)(t,\cdot)\|_{L^2(\mathbb{R}^N,\rho)}$ for any $t\in\mathbb{R}_\tau$, we can write the inequality in (\ref{Des_F}) of the following form
\[
\varphi(t)\leqslant M_1\| h_\epsilon- h_0\|_{L^\infty(\mathbb{R})}e^t+M_1\|J\|_{L^\infty(\mathbb{R}^N)}\rho_1^{-1}\int_\tau^t\varphi(s)ds,
\]
and using Gronwall's inequality
\[
\|(u^\epsilon -u^0)(t,\cdot)\|_{L^2(\mathbb{R}^N,\rho)}\leqslant M_1\| h_\epsilon- h_0\|_{L^\infty(\mathbb{R})}e^{M_1\|J\|_{L^\infty(\mathbb{R}^N)}\rho_1^{-1}(t-\tau)}
\]
where $M_1=2^{(p+1)/p} \ell_g\beta>0$ and $\rho_1=\inf\{\rho(\xi);\ |\xi|\leqslant 1 \}>0$, concluding the proof. $\blacksquare$

\begin{theorem}
The pullback attractor $\{\mathcal{A}_\epsilon(t);\ t\in\mathbb{R}\}$ is upper semicontinuous in $\epsilon=0$.
\end{theorem}
\noindent  $\proof$ Let $\tau\in\mathbb{R}$ such that $\mathrm{dist}(S_0(t,\tau)B(0;a),\mathcal{A}_0(t))<\delta/2$, where $\bigcup_{s\in\mathbb{R}}\mathcal{A}(s)\subset B(0;a)$ for all $\delta>0$. By Theorem \ref{The01}
$$
\|S_\epsilon(t,\tau)u_0-S_0(t,\tau)u_0\|_{L^2(\mathbb{R}^N,\rho)}\leqslant   M_1\| h_\epsilon- h_0\|_{L^\infty(\mathbb{R})}e^{M_1\|J\|_{L^\infty(\mathbb{R}^N)}\rho_1^{-1}(t-\tau)}\to0
$$
as $\epsilon\to0^+$ in compact subsets of $\mathbb{R}$ uniformly for $u_0$ in bounded subsets of $L^2(\mathbb{R}^N,\rho)$. Hence, there exists $\epsilon_0>0$ such that
\[
\sup_{a_\epsilon\in\mathcal{A}_\epsilon(\tau)}\|S_\epsilon(t,\tau)a_\epsilon-S_0(t,\tau)a_\epsilon\|_{L^2(\mathbb{R}^N,\rho)}<\delta/2,\ \forall\ \epsilon\in[0,\epsilon_0].
\]

Then
\begin{eqnarray*}
 &&\mathrm{dist}(\mathcal{A}_\epsilon(t),\mathcal{A}_0(t))\\
&\leqslant& \mathrm{dist}(S_\epsilon(t,\tau)\mathcal{A}_\epsilon(\tau),S_0(t,\tau)\mathcal{A}_\epsilon(\tau))+\mathrm{dist}(S_0(t,\tau)\mathcal{A}_\epsilon(\tau),S_0(t,\tau)\mathcal{A}_0(\tau))\\
&=&\sup_{a_\epsilon\in\mathcal{A}_\epsilon(\tau)} \mathrm{dist}(S_\epsilon(t,\tau)a_\epsilon,S_0(t,\tau)a_\epsilon)+\mathrm{dist}(S_0(t,\tau)\mathcal{A}_\epsilon(\tau),\mathcal{A}_0(t))<\frac{\delta}{2}+\frac{\delta}{2}=\delta
\end{eqnarray*}
and the upper semicontinuity is proved. $\blacksquare$

\begin{remark}
Assuming that $g$ is a function of class $\mathcal{C}^k$, for any integer $k\geqslant0$,  and its derivatives up to order $k$ are bounded we can use the estimates in $L^2(\mathbb{R}^N,\rho)$ have been obtained for solution in the pullback attractor and a bootstrap argument to obtain $\mathcal{C}^k$ estimates.
\end{remark}

\begin{remark}
As well as in the autonomous case, $h(t)\equiv h$, it follows from Theorem \ref{Theorem6} that the pullback attractor is contained in the space of bounded continuous function.
\end{remark}

\end{document}